\documentclass{amsart}

\newtheorem{theorem}{Theorem}
\newcommand{\bt}{\begin{theorem}}
\newcommand{\et}{\end{theorem}}
\newtheorem{lemma}{Lemma}
\newcommand{\bl}{\begin{lemma}}
\newcommand{\el}{\end{lemma}}
\newtheorem{corollary}{Corollary}
\newcommand{\bc}{\begin{corollary}}
\newcommand{\ec}{\end{corollary}}
\newtheorem{problem}{Problem}
\newcommand{\bprob}{\begin{problem}}
\newcommand{\eprob}{\end{problem}}

\newtheorem{example}{Example}
\newcommand{\bex}{\begin{example}}
\newcommand{\eex}{\end{example}}

\newcommand{\beq}{\begin{equation}}
\newcommand{\eeq}{\end{equation}}
\newcommand{\benum}{\begin{enumerate}}
\newcommand{\eenum}{\end{enumerate}}
\newcommand{\N}{\ensuremath{ \mathbf N }}

\newcommand{\mcr}{\ensuremath{ \mathcal R}}

\newcommand{\bsmallmat}{\left(\begin{smallmatrix}}
\newcommand{\esmallmat}{\end{smallmatrix}\right)}

\newcommand{\bmat}{\left(\begin{matrix}}
\newcommand{\emat}{\end{matrix}\right)}

\DeclareMathOperator{\qqand}{\qquad\text{and}\qquad}

\usepackage[all]{xy}

\usepackage{amsfonts}
\usepackage{amsmath,amssymb,amsthm}
\usepackage{color}

\title[Iterated Syracuse function]{Permutation Patterns of the  Iterated Syracuse Function}
\author{Melvyn B. Nathanson}
\address{Department of Mathematics\\Lehman College (CUNY)\\Bronx, NY 10468} 
\email{melvyn.nathanson@lehman.cuny.edu}



\date{\today}

\begin{document}

\begin{abstract}
Let $\Omega$ be  the set of odd positive integers and let $S:\Omega \rightarrow \Omega$ 
be the  Syracuse function. 
It is proved that, for every permutation $\sigma$ of $(1,2,3)$, 
the set of triples of the form $(m,S(m),S^2(m))$ 
with permutation pattern $\sigma$ has positive density, and these densities are computed.  
However, there exist permutations $\tau$ of $(1,2,3,4)$ such that no quadruple 
$(m,S(m), S^2(m), S^3(m))$ has permutation pattern $\tau$.  
This implies the nonexistence of certain permutation patterns of $n$-tuples 
$(m,S(m),\ldots, S^{n-1}(m))$ for all $n \geq 4$. 
\end{abstract}

\maketitle

\section{Permutation patterns}

Let $\Omega$ be  the set of odd positive integers. 
The \emph{Syracuse function}\index{Syracuse function} is the arithmetic function 
$S:\Omega \rightarrow \Omega$  
defined by 
\[
S(m) = \frac{3m + 1}{2^e}
\]
where $e$ is the largest integer such that $2^e$ divides  $3m+1$. 
Equivalently, $S(m)$ is the \emph{odd part}, that is, 
the largest odd divisor, of the even integer $3m+1$. 
Note that $S(m) = 1$ if and only if $m = \sum_{i=0}^k 4^i  = (4^{k+1}-1)/3$ 
for some nonnegative integer $k$.  
The notorious \emph{Collatz conjecture}\index{Collatz conjecture} 
asserts that for every positive integer $m$ 
there exists an integer $r$ such that 
$S^r(m) = 1$  (Lagarias~\cite{laga85}, Tao~\cite{tao22}, Wikipedia~\cite{wiki}).  
By supercomputer calculation, Barina~\cite{bari21} has verified  the conjecture  
for all $m < 2^{68}$.  

An arithmetic function is any function whose domain is a 
nonempty subset $\Omega$ of the set positive integers.   
Let $S: \Omega \rightarrow \Omega$ be an arithmetic function 
and, for $j \in \N$, let $S^j:\Omega \rightarrow \Omega$  be the $j$th iterate of $S$.   
Let $V = (v_i)_{i=1}^n$ be a finite sequence of positive integers.  
We say that an integer $m$ in $\Omega$ has \emph{increasing-decreasing pattern}\index{pattern} 
$V$ with respect to $S$ if 
\[
m < S(m) < S^2(m) < \cdots < S^{v_1}(m) 
\]
\[
S^{v_1}(m) > S^{v_1+1}(m) > \cdots > S^{v_1+v_2}(m)
\]
\[
S^{v_1+v_2}(m) < S^{v_1+v_2+1}(m) < \cdots < S^{v_1+v_2+v_3}(m)
\]
and, in general, if $i$ is odd, then 
\beq                                                                 \label{Syracuse:increasing} 
S^{v_1+ \cdots + v_{i-1}}(m) < S^{v_1+ \cdots + v_{i-1}+1}(m) < \cdots < S^{v_1+ \cdots + v_{i-1}+v_{i}}(m)
\eeq
and if $i$ is even, then 
\beq                                                                 \label{Syracuse:decreasing} 
S^{v_1+ \cdots + v_{i-1}}(m) > S^{v_1+ \cdots  +v_{i-1}+1}(m) > \cdots > S^{v_1+ \cdots  +v_{i-1}+v_i}(m).
\eeq
The arithmetic function $S$ is \emph{wildly increasing-decreasing} if,  
for every finite sequence $V$ of positive integers, there exists an integer $m \in \Omega$ 
such that $m$ has increasing-decreasing pattern $V$ with respect to $S$.  

Nathanson~\cite{nath24z} proved that the Syracuse function is 
wildly increasing-decreasing.  In this paper we consider more subtle 
variations in successive iterates of the Syracuse function.

Let $\Sigma_n$ be the group of permutations of $\{1,2,3,\ldots, n\}$. 
Let $X = (x_1, x_2,\ldots, x_n)$ be an $n$-tuple of  distinct real numbers. 
We rearrange the coordinates of $X$ to obtain an $n$-tuple
$ (y_1, y_2, \ldots, y_n)$ such that 
\[
y_1 < y_2< \cdots < y_n.
\]
There is a unique permutation $\sigma \in \Sigma_n$ such that 
\[
 (x_1, x_2,\ldots, x_n) =  (y_{\sigma(1)}, y_{\sigma(2)}, \ldots, y_{\sigma(n)}). 
\]
We call $\sigma$ the  \emph{permutation pattern}\index{permutation pattern} of the $n$-tuple $X$
and denote it by $(\sigma(1), \sigma(2), \ldots, \sigma(n))$.  
In standard form, this is the permutation 
$\bsmallmat 1 & 2 & \cdots & n \\ \sigma(a) & \sigma(2) & \cdots & \sigma(n) \esmallmat$. 

For example, if 
\[
X = (x_1, x_2, x_3, x_4) = (7, 13, 18, 11)
\]
 then 
\[
7 < 11 < 13 < 18. 
\]
We obtain  
\[
(y_1, y_2, y_3, y_4) = ( 7, 11, 13, 18).
\]
and 
\[
(x_1, x_2, x_3, x_4) = (y_1, y_3, y_4, y_2). 
\]
The permutation pattern of the quadruple $(7, 13, 18, 11)$ 
is $\sigma =  ( 1 ,3,4,2)$.

It is an open problem to determine, for every positive integer $n$,  
the possible permutation patterns of the initial segments  
$\left( m, S(m), S^2(m),\ldots, S^{n-1}(m) \right)$ of the iterated Syracuse function 
for  integers $m \in \Omega$ such that $S^i(m) \neq S^j(m)$ for $0 \leq i < j \leq n-1$. 
For every permutation $\sigma \in \Sigma_n$, let $\Gamma_{\sigma}(M)$ count the number of odd 
positive integers $m \leq M$ such that the $n$-tuple 
$\left( m, S(m), S^2(m),\ldots, S^{n-1}(m) \right)$ 
has distinct coordinates and has permutation pattern $\sigma$.  
The \emph{permutation density}\index{permutation density} of $\sigma \in \Sigma_n$  is
\[
d_n(\sigma) = \lim_{M\rightarrow \infty} \frac{\Gamma_{\sigma}(M)}{M/2}
\]
(if the limit exists).

In this paper we prove that every permutation $\sigma \in \Sigma_3$ 
occurs with positive density.  
We also prove that there are permutations $\tau \in \Sigma_4$ 
such that  $\tau$ not only has zero permutation density, 
but there exists no positive integer $m$ with permutation pattern $\tau$.  
If no $n$-tuple  $\left( m, S(m), S^2(m),\ldots, S^{n-1}(m) \right)$ has permutation 
pattern $(a_1,\ldots, a_n) \in \Sigma_n$, 
then no $(n+1)$-tuple  $\left( m, S(m), S^2(m),\ldots, S^{n-1}(m) , S^n(m)\right)$ 
has permutation pattern $(a_1,\ldots, a_n, n+1) \in \Sigma_{n+1}$. 
It follows that, for all $n \geq 4$, 
there exist permutations $\tau$  in $\Sigma_n$ such that 
no $n$-tuple  $\left( m, S(m), S^2(m),\ldots, S^{n-1}(m) \right)$ 
has permutation pattern $\tau$.

\bt                      \label{Syracuse:theorem:triples}
The following table gives the density of  permutation patterns of triples $\left( m, S(m), S^2(m) \right)$.
\begin{center}
\begin{tabular}{ c|c}
permutation pattern $\sigma \in \Sigma_3$ & permutation density $d_3(\sigma) $ \\     \hline 
$(1,2,3)$ & $1/4$ \\
$(1,3,2)$ & $1/8$ \\
$(2,1,3)$ & $1/8$ \\
 $(2,3,1)$ & $1/8$  \\
$(3,1,2)$ & $1/8$  \\ 
$(3,2,1)$ & $1/4$ \\
 \end{tabular}
\end{center}
\et

The proof of this result will follow immediately from Theorems~\ref{Syracuse:theorem:triples-not-5} 
and~\ref{Syracuse:theorem:triples-5} below. 

Note: 
In the study of the Collatz conjecture,  instead of the Syracuse function, 
investigators often use the Collatz functions 
\[
C(m) = \begin{cases}
3m+1 & \text{if $m$ is odd} \\
\frac{m}{2} & \text{if $m$ is even}
\end{cases}
\]
and
\[
C_1(m) = \begin{cases}
\frac{3m+1}{2} & \text{if $m$ is odd} \\
\frac{m}{2} &  \text{if $m$ is even.}
\end{cases}
\]
The equivalent Collatz conjectures state that for every odd integer $m$ there exist  
positive integers $n$ and $n_1$ such that $C^n(m) = 1$ and $C_1^{n_1}(m) = 1$. 
Simons and  de Weger~\cite{simo-dewe05} 
and Simons~\cite{simo07,simo08a,simo08b}
have studied a different kind of increasing-decreasing behavior for the Collatz function. 
 
\section{Permutation patterns  for pairs}
We begin with a simple but important calculation of permutation densities for pairs. 

\bt                    \label{Syracuse:theorem:pairs}
Let $m \in \Omega$. 
If $m \equiv 3 \pmod{4}$, then $(m,S(m))$ has permutation pattern $(1,2)$.  
If $m \equiv 1 \pmod{4}$ and $m > 1$, then $(m,S(m))$ has permutation pattern $(2,1)$.   
\begin{center}
\begin{tabular}{ c|c|c}
congruence class of $m$ &  permutation pattern $\sigma \in \Sigma_2$ 
& permutation density $d_3(\sigma)$ \\
     \hline 
$3\pmod{4}$  & $(1,2)$ & $1/2$ \\
$1\pmod{4}$  & $(2,1)$ & $1/2$ \\
 \end{tabular}
\end{center}
\et

\begin{proof}
If $m = 3+4x$ for some nonnegative integer $x$, then 
$3m+1=10+12x$ and 
\[
S(m) = \frac {10+12x}{2} = 5+6x >  3+4x = m.
\]
Thus, $(m,S(m))$ has permutation pattern $(1,2)$.

If $m = 1+4x$ for some positive integer $x$, then 
$3m+1=4+12x$ and there is an integer $e \geq 2$ such that 
\[
S(m) = \frac{4+12x}{2^e} = \frac{1+3x}{2^{e-2}} \leq 1+3x < 1+4x = m
\]
and $(m,S(m))$ has permutation pattern $(2,1)$.
This completes the proof. 
\end{proof}

\bc
Let $A$ be the set of odd positive integers $m$ for which the permutation pattern 
of the triple $(m,S(m), S^2(m))$ is 
\[
 (1,2,3), \qquad (1,3,2), \quad \text{or} \quad (2,3,1).
\]
Let $B$ be the set of odd positive integers $m$ for which the permutation pattern 
of the triple $(m,S(m), S^2(m))$ is 
\[
 (2,1,3), \qquad (3,1,2), \quad \text{or} \quad (3,2,1).
\]
The set $A$ has density $1/2$ and the set $B$ has density $1/2$. 
\ec

\begin{proof}
The pair $(m,S(m))$ has permutation pattern $(1,2)$ if and only if 
the triple $(m,S(m),S^2(m))$ has permutation pattern 
$ (1,2,3)$, $(1,3,2)$, or $(2,3,1)$. 
Similarly, the  pair $(m,S(m))$ has permutation pattern $(2,1)$ if and only if 
the triple $(m,S(m),S^2(m))$ has permutation pattern 
$ (2,1,3)$,  $(3,1,2)$, or $(3,2,1)$. 
This completes the proof. 
\end{proof}

\section{Permutation patterns for triples for $m \not\equiv 5 \pmod{8}$}
The calculation of permutation densities for triples $\left( m, S(m), S^2(m) \right)$ 
is divided into two cases.
In the first case we consider odd positive integers $m \not\equiv 5 \pmod{8}$ 
and in the second case we consider odd positive integers $m \equiv 5 \pmod{8}$.

\bt          \label{Syracuse:theorem:triples-not-5}
Every odd positive integer $m$ such that $m \not\equiv 5 \pmod{8}$ 
belongs to exactly one of the five congruence classes in the table below.  
Each congruence class uniquely determines the 
permutation pattern of the triple $(m, S(m), S^2(m))$  for all integers $m > 1$ 
in the congruence class. 

\vspace{0.4cm}

\begin{center}
\begin{tabular}{ c|c|c} 
congruence class of $m$ &  permutation pattern $\sigma \in \Sigma_3$ 
& permutation density  $d_3(\sigma)$ \\     \hline 
$7\pmod{8} $   & $(1,2,3)$ & $1/4$ \\  
$9 \pmod{16}$ &  $(2,1,3)$ & $1/8$ \\ 
$11 \pmod{16}$ &   $(1,3,2)$  & $1/8$\\ 
$3 \pmod{16}$ &   $(2,3,1)$ & $1/8$ \\ 
$1 \pmod{16}$ &  $(3,2,1)$ & $1/8$ 
\end{tabular}
\end{center}
\et

\begin{proof}
The following diagram shows the five congruence classes (in black boxes)
that partition the  integers $m \not\equiv 5 \pmod{8}$
and the ``missing'' congruence class $5 \pmod{8}$ (in the red box).

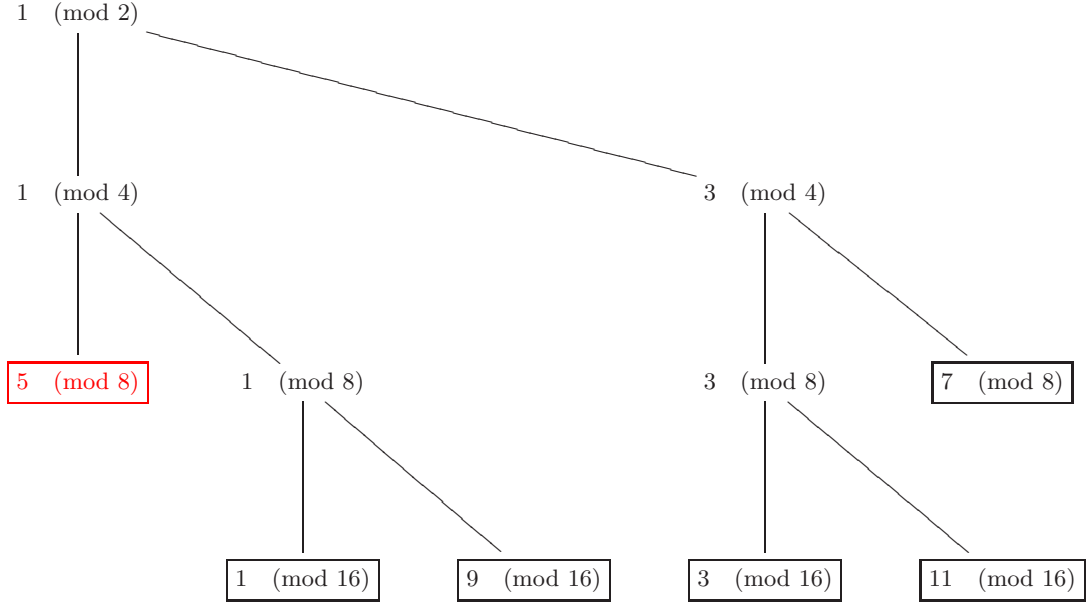
\small{
\[
\xymatrix{
 1 \pmod {2} \ar@{-}[dd]   \ar@{-}[ddrrr]   & & &  &   \\
 & & &  &    \\
1\pmod{4}   \ar@{-}[dd]   \ar@{-}[ddr]   & &  & 3 \pmod{4}  \ar@{-}[dd] \ar@{-}[ddr] & \\
& & & &   \\
\color{red}\boxed{5 \pmod{8} } &  1 \pmod{8}  \ar@{-}[dd]   \ar@{-}[ddr]   & &  3 \pmod{8}   \ar@{-}[dd]   \ar@{-}[ddr]  & \boxed{7 \pmod{8}} \\
 & & & & &   \\
 &\boxed{1 \pmod{16}} &  \boxed{9 \pmod{16}}  &   \boxed{3\pmod{16}}  & \boxed{11 \pmod{16}} \\
}
\]
}

\normalsize{ 

If $m \equiv 7 \pmod{8}$, 
then for some  nonnegative integer $x$ we have 
\begin{align*}
m & = 7 + 8x \\
S(m) & = \frac{22+24x}{2} = 11 + 12 x \\
S^2(m) & = \frac{34+36x}{2} =  17 + 18x 
\end{align*}
and 
\[
7 + 8x < 11 + 12 x < 17 + 18x.
\]
Therefore, $(m, S(m), S^2(m))$ has permutation pattern $(1,2,3)$. 

If $m \equiv 9 \pmod{16}$, 
then for some  nonnegative integer $x$ we have 
\begin{align*}
m & = 9 + 16x \\
S(m) & = \frac{28+48x}{4} = 7 + 12 x \\
S^2(m) & = \frac{22+36x}{2} =  11 + 18x 
\end{align*}
and 
\[
7 + 12x < 9 + 16 x < 11 + 18x.
\]
Therefore, $(m, S(m), S^2(m))$ has permutation pattern $(2,1,3)$.

If $m \equiv 11 \pmod{16}$, 
then for some  nonnegative integer $x$ we have 
\begin{align*}
m & = 11 + 16x \\
S(m) & = \frac{34+48x}{2} = 17 + 24 x \\
S^2(m) & = \frac{52+72x}{4} =  13 + 18x 
\end{align*}
and 
\[
 11 + 16 x < 13 + 18x < 17 + 24x.
\]
Therefore, $(m, S(m), S^2(m))$ has permutation pattern $(1,3,2)$.

If  $m \equiv 3 \pmod{16}$, 
then for some positive integer $x$ and nonnegative integer $e$ we have 
\begin{align*}
m & = 3 + 16x \\
S(m) & = \frac{10+48x}{2} = 5 + 24 x \\
S^2(m) & = \frac{16+72x}{8} = \frac{2 + 9x}{2^e }.
\end{align*}
and 
\[
\frac{2 + 9x}{2^e } \leq 2+9x  < 3 + 16x < 5 + 24 x.
\]
Therefore, $(m, S(m), S^2(m))$ has permutation pattern $(2,3,1)$.

If $m > 1$ and $m \equiv 1 \pmod{16}$, 
then for some positive integer $x$ and nonnegative integer $e$ we have 
\begin{align*}
m & = 1 + 16x \\
S(m) & = \frac{4+48x}{4} = 1 + 12 x \\
S^2(m) & = \frac{4+36x}{4} = \frac{1 + 9x}{2^e }.
\end{align*}
and 
\[
\frac{1 + 9x}{2^e } \leq 1 + 9x < 1 + 12 x < 1 + 16x.
\]
Therefore,  $(m, S(m), S^2(m))$ has permutation pattern $(3,2,1)$.
This completes the proof. 
}
\end{proof}


\section{Permutation patterns  for triples  for $m \equiv 5 \pmod{8}$}

We shall prove that half of the odd positive integers congruent to $5 \pmod{8}$ 
have the 3-term permutation pattern 
$(3,2,1)$ and half have the 3-term permutation pattern $(3,1,2)$.

In the following proofs, $e$ denotes a nonnegative integer.

\bl                                           \label{Syracuse:lemma:1}
Let $k \geq 1$.  If $m$ is a positive integer and 
\[
m \equiv  \sum_{i=0}^{k} 4^i  +  2\cdot 4^k  \pmod{ 2 \cdot 4^{k+1}} 
\]
then $m \equiv 5 \pmod{8}$ and the permutation pattern of $(m,S(m), S^2(m))$ is $(3,2,1)$.
\el

\begin{proof}
For some nonnegative integer $x$ we have 
\[ 
m  = \frac{4^{k+1}-1}{3} + 2\cdot 4^k  + 2 \cdot 4^{k+1}  x  \geq 13 + 32x 
\] 
 and  
\[
3m+1 = 4^{k+1} + 3\cdot 2 \cdot 4^{k} + 3 \cdot 2\cdot 4^{k+1} x. 
\]
It follows that 
\[
S(m) =  \frac{3m+1}{2 \cdot 4^k } = 5 + 12x < 13+  32x \leq m.
\]
We have 
\[
3S(m)+1 = 16+36x = 4(4+9x) 
\]
and so 
\[
S^2(m) = \frac{4 +9x}{2^e} \leq  4 +9x < 5 + 12x = S(m). 
\]
Therefore, 
\[
S^2(m) < S(m) < m
\]
and the permutation pattern of $(m,S(m), S^2(m))$ is $(3,2,1)$. 
This completes the proof. 
\end{proof}

\bl                                             \label{Syracuse:lemma:2}
Let $k \geq 1$.  If $m$ is a positive integer and 
\[
m \equiv \sum_{i=0}^{k} 4^i   \pmod{4^{k+2}}
\]
then  $m \equiv 5 \pmod{8}$ and either $m = \sum_{i=0}^{k} 4^i$ 
and $S(m) = S^2(m) = 1$ 
or the permutation pattern of $(m,S(m), S^2(m))$ is $(3,2,1)$. 
\el

\begin{proof}
For some nonnegative integer $x$ we have 
\[
m  = \frac{4^{k+1}-1}{3}+ 4^{k+2} x  \geq 5 + 64x 
\] 
and  
\[
3m+1 = 4^{k+1} + 3\cdot 4^{k+2} x.  
\]
If $x = 0$, then $3m+1 = 4^{k+1}$ and $S(m) = S^2(m) = 1$.

If $x \geq 1$, then 
\[
S(m) = \frac{3m+1}{4^{k+1}} = 1 + 12x < 5 + 64 x \leq m 
\]
and  
\[
S^2(m) = \frac{1+9x}{2^e} \leq 1+9x < 1 + 12x = S(m). 
\]
Therefore, 
\[
S^2(m) < S(m) < m
\]
and the  permutation pattern of $(m,S(m), S^2(m))$ is $(3,2,1)$. 
This completes the proof. 
\end{proof}

\bl                                                     \label{Syracuse:lemma:3}
Let $k \geq 1$.   If $m$ is a positive integer and 
\[
m \equiv  \sum_{i=0}^{k+1} 4^i   + 2\cdot 4^{k}   \pmod{ 2\cdot 4^{k+1} }
\]
then  $m \equiv 5 \pmod{8}$ and  the permutation pattern of $(m,S(m), S^2(m))$ is $(3,1,2)$.
\el

\begin{proof}
For some nonnegative integer $x$ we have 
\begin{align*}
m & = \frac{4^{k+2}-1}{3}  + 2\cdot  4^{k}  + 2 \cdot 4^{k+1}  x \\
& \geq 29 + 32x 
\end{align*}  
and  
\[
3m+1 = 4^{k+2} + 3\cdot 2\cdot  4^{k}  + 3 \cdot 2 \cdot 4^{k+1} x.
\]
Therefore,
\[
S(m) = \frac{3m+1}{ 2\cdot 4^{k}} = 11 + 12x 
\]
and 
\[
S^2(m) = 17+18x.
\] 
We have 
\[
S(m) = 11+12x < S^2(m) = 17+18x < 29+32x \leq m
\]
and the permutation pattern of $(m,S(m), S^2(m))$ is $(3,1,2)$. 
This completes the proof. 
\end{proof}

\bl                                                          \label{Syracuse:lemma:4}
Let $k \geq 1$.  If $m$ is a positive integer and 
\[
m \equiv \sum_{i=0}^{k+1} 4^i +  4^{k+1}  \pmod{4^{k+2}}
\]
then $m \equiv 5 \pmod{8}$ and  the permutation pattern of $(m,S(m), S^2(m))$ is $(3,1,2)$.
\el

\begin{proof}
For some nonnegative integer $x$ we have 
\begin{align*}
m & = \frac{4^{k+2}-1}{3}  + 4^{k+1}+ 4^{k+2} x \\
& \geq 37 + 64x 
\end{align*}  
and  
\[
3m+1 = 4^{k+2} + 3\cdot 4^{k+1} + 3 \cdot 4^{k+2}  x 
\]
It follows that 
\[
S(m) =  \frac{3m+1}{4^{k+1} } = 7 + 12x 
\]
and 
\[
S^2(m) = 11+18x. 
\]
Therefore, 
\[
S(m)  = 7 + 12 < S^2(m) = 11+18x < 37 + 64x \leq m
\]
and the permutation pattern of $(m,S(m), S^2(m))$ is $(3,1,2)$.
This completes the proof. 
\end{proof}

\emph{Notation.} 
Let $F(M)$ be a positive function of $M$.  
We denote by  $o(F(M))$ a function $G(M)$  such that 
 $\lim_{M\rightarrow \infty} \frac{G(M)}{F(M)} = 0$. 
Thus, $o(1)$ denotes a function $G(M)$ such that $\lim_{M\rightarrow \infty} G(M)  = 0$. 
Note that $-o(1) = o(1)$ and $o(1) + o(1) = o(1)$.  
Also, $o(X(M)) = X(M)o(1)$. 

The \emph{counting function} of a set $A$ of positive integers is 
\[
A(M) = \sum_{\substack{ a \in A \\ a \leq M}} 1. 
\]
A subset $A$ of a set $X$ of positive integers has density $\alpha$ 
with respect to $X$ if the limit
\[
\lim_{M \rightarrow \infty} \frac{A(M)}{X(M)} 
\]
exists and equals $\alpha$.  
Equivalently, $A$ has density $\alpha$ with respect to $X$, denoted $d_X(A) = \alpha$, 
 if  
\[
A(M) = \alpha X(M) + o(X(M)). 
\] 
If $d_X(A) = 0$, then $A(M) = o(X(M))$.


\bl                                                        \label{Syracuse:lemma:densityW}
Let $X$ be a set of positive integers, let $W$ be a subset of $X$ with 
density $d_X(W) = \omega > 0$, and let 
\[
W = R_0 \cup R_1 \cup \cdots \cup  R_t 
\]
be a partition of $W$. 
Let $\alpha_1,\ldots, \alpha_t$ be positive real numbers such that 
$\sum_{i=1}^t \alpha_i = \omega$.  
If $d_X(R_0) = 0$ and if, for all $\varepsilon$ with 
\[
0 < \varepsilon < \min(\alpha_1,\ldots, \alpha_t) 
\]
there is a subset $R_{i,\varepsilon}$ of $R_i$ such that $d_X(R_{i,\varepsilon}) = \alpha_i- \varepsilon$, then $d_X(R_j) = \alpha_j$ for all $j \in \{1,\ldots, t\}$.
\el

\begin{proof}
Let $M \geq 1$.  Because the sets $R_0, R_1, \ldots, R_t$ partition the set $W$, 
we have the counting function equation  
\[
W(M) = R_0(M) + R_1(M) + \cdots + R_t(M). 
\]
The density condition $d_X(R_0) = 0$ implies 
\[
R_0(M) =  o(X(M)). 
\]
Let $0 < \varepsilon < \min(\alpha_1,\ldots, \alpha_t)$.  
For all $i \in \{1,\ldots, t\}$, 
the subset condition $R_{i,\varepsilon} \subseteq R_i$ implies 
\[
0 \leq R_{i,\varepsilon}(M) \leq R_i(M). 
\]
The density condition 
\[
d_X(R_{i,\varepsilon}) = \lim_{M\rightarrow \infty} \frac{R_{i,\varepsilon}(M)}{X(M)} 
= \alpha_i - \varepsilon
\]
implies 
\[
R_{i,\varepsilon}(M) = ( \alpha_i - \varepsilon + o(1))X(M). 
\]
Also, $d_X(W) = \omega > 0$ implies 
\[
W(M) = (\omega + o(1))X(M).
\]
For all $j \in \{1,\ldots, t\}$, we have 
\begin{align*}
\left( \alpha_j - \varepsilon  + o(1) \right) X(M)  
& =  R_{j,\varepsilon}(M) \\ 
&  \leq R_j(M) \\
& = W(M) - \sum_{\substack{i=0 \\ i \neq j }}^t R_i(M) \\
& \leq  W(M) -  R_0(M) - \sum_{\substack{i=1 \\ i \neq j }}^t  R_{i,\varepsilon}(M) \\ 
& =  (\omega + o(1))X(M)- o(1)X(M)   
- \sum_{\substack{i=1 \\ i \neq j}}^t  \left( ( \alpha_i - \varepsilon +o(1))X(M) \right) \\ 
& =  \omega X(M) - \sum_{\substack{i=1 \\ i \neq j}}^t  ( \alpha_i - \varepsilon)X(M) + o(1)X(M) \\ 
& =  X(M) \left( \omega - \sum_{\substack{i=1 \\ i \neq j}}^t  \alpha_i + (t-1) \varepsilon + o(1) \right) \\
& =  X(M) \left(  \alpha_j + (t-1) \varepsilon + o(1) \right).  
\end{align*}
Therefore, 
\[
\alpha_j - \varepsilon + o(1)  \leq \frac{R_j(M) }{X(M)} 
\leq  \alpha_j + (t-1) \varepsilon  + o(1)
\]
for all $\varepsilon > 0$ and so 
\[
\alpha_j + o(1)  \leq \frac{R_j(M) }{X(M)} 
\leq  \alpha_j  + o(1).
\]
Thus, 
\[
d_X(R_j) = \lim_{M \rightarrow \infty} \frac{R_j(X(M)) }{X(M)} =  \alpha_j 
\]
for all  $j \in \{1,\ldots, t\}$. 
This completes the proof. 
\end{proof}


\bt                                                  \label{Syracuse:theorem:triples-5}
Let $m$ be a  positive integer such that $m \equiv 5 \pmod{8}$. 
If $m = \sum_{i=0}^k 4^i$ for some $k \geq 1$, then $\left(m, S(m), S^2(m)\right) = (m,1,1)$.
If $m \neq \sum_{i=0}^k 4^i$ for some $k \geq 1$, 
then the triple $\left(m, S(m), S^2(m)\right)$ 
has permutation pattern $(3,2,1)$ or $(3,1,2)$.
In the congruence class  $m \equiv 5 \pmod{8}$, 
the permutation  densities are as follows:  
\begin{center}
\begin{tabular}{ c|c}
permutation pattern $\sigma \in \Sigma_3$ & permutation density  $d_3(\sigma)$ \\     \hline 
$(3,2,1)$ & $1/8$ \\
$(3, 1,2)$ & $1/8$ \\
 \end{tabular}
\end{center}

\et

\begin{proof} 
For every positive integer $k$, let 
\[
r_k = \sum_{i=0}^k 4^i = \frac{4^{k+1}-1}{3}
\]
and let 
\[
\mcr_0  =  \{r_k:k=1,2,3\ldots\} = \{5,21, 85, 341, 1365,\ldots \}.
\]
The set  $\mcr_0$ has density zero.  
We have $S(m) = 1$ if and only if $m=1$ or $m \in \mcr_0$.  
If $m \in \mcr_0$, then $m \equiv 5 \pmod{8}$.  

Let 
\[
\mcr_1 = \left\{m \equiv 5 \pmod{8}: \text{$\left(m, S(m), S^2(m)\right)$ 
has permutation pattern $(3,1,2)$} \right\}
\]
and 
\[
\mcr_2 = \left\{m \equiv 5 \pmod{8}: \text{$\left(m, S(m), S^2(m)\right)$ 
has permutation pattern $(3,2,1)$} \right\}.
\]
The sets $\mcr_0, \mcr_1$, and $\mcr_2$ are pairwise disjoint.

Let 
\[
5\pmod{8} \setminus \left( \mcr_0 \cup r_{k+1} \pmod{2\cdot 4^{k+1}} \right) 
\]
\[
 \left( 5\pmod{8} \setminus \mcr_0 \right)  
  \setminus \left(  r_{k+1} \pmod{2\cdot 4^{k+1}} \right) 
\]
denote the set of positive integers $m \notin \mcr_0$ such that 
$m \equiv 5\pmod{8} $ but $m \not\equiv r_{k+1} \pmod{2\cdot 4^{k+1}}$. 
We shall prove by induction on $k$  that  this set 
is partitioned into two sets, one with permutation pattern $(3,2,1)$ 
and the other with permutation pattern $(3,1,2)$, and that each of these sets 
has permutation density
\[
\frac{1}{8}\left(1 -  \frac{1}{4^k} \right).
 \]

We begin with the cases $k=1$ and $k=2$.

Here is a picture of a partition of the congruence class $5\pmod{8}$ 
into disjoint unions of congruence classes.

\newpage

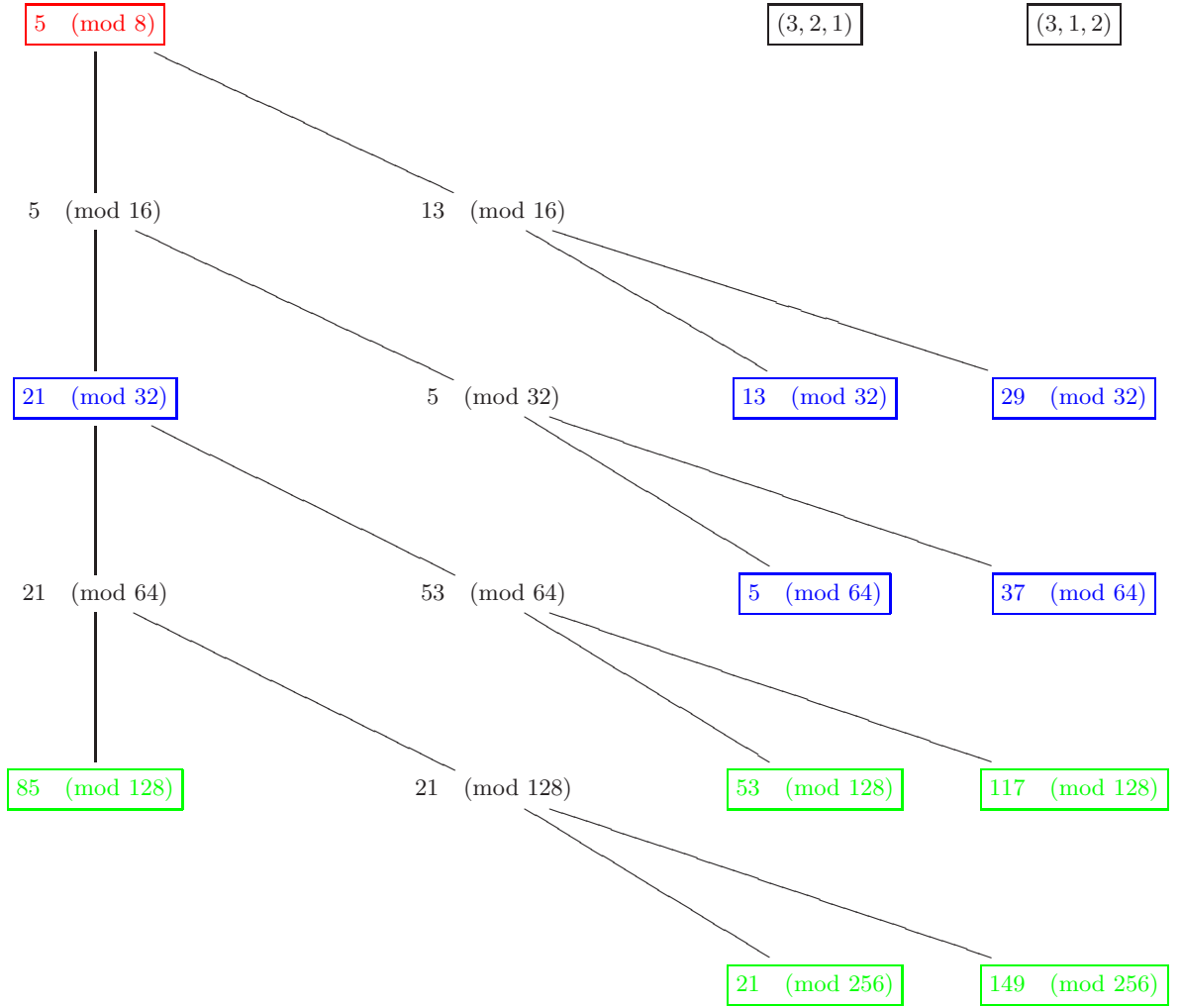
\small{
\[
\xymatrix{
\color{red}\boxed{ 5 \pmod {8} } \ar@{-}[dd]   \ar@{-}[ddrrr]  & & & & &  \boxed{{(3,2,1)}} &  \boxed{{(3,1,2)}}  \\
 & & & & &  &   \\
5\pmod{16}   \ar@{-}[dd]   \ar@{-}[ddrrr]  &  & &  13 \pmod{16}  \ar@{-}[ddrr]   \ar@{-}[ddrrr] & & & \\
& & & & & &   \\
 \color{blue}\boxed{ 21\pmod{32} }\ar@{-}[dd]   \ar@{-}[ddrrr] & & & 5 \pmod{32}  \ar@{-}[ddrr]   \ar@{-}[ddrrr]   & &  
  \color{blue}\boxed{13 \pmod{32}}  & \color{blue}\boxed{29 \pmod{32}} \\
  & & & & & &   \\
   21\pmod{64} \ar@{-}[dd]   \ar@{-}[ddrrr] & & & 53 \pmod{64}   \ar@{-}[ddrr]   \ar@{-}[ddrrr]   & &  \color{blue}\boxed{5\pmod{64}}  & \color{blue}\boxed{37 \pmod{64}} \\
  & & & & & &   \\
\color{green}\boxed{85 \pmod{128} } & & & 21 \pmod{128}   \ar@{-}[ddrr]   \ar@{-}[ddrrr]    & &  \color{green}\boxed{53 \pmod{128}}  & \color{green}\boxed{117 \pmod{128} } \\
  & & & & & &   \\
& & &  & &  
  \color{green}\boxed{21 \pmod{256}}  & \color{green}\boxed{149 \pmod{256}} \\
  & & & & & &   \\ 
}
\]
}

\normalsize{ }

The congruence class $5 \pmod{8} =  r_1 \pmod{8}$ is the disjoint union of the following 
 five congruence classes (in the blue boxes):
\begin{align*}
13\pmod{32} & \\
5\pmod{64} &  \\
29\pmod{32} &  \\
37 \pmod{64}& \\
21\pmod{32} & = r_2 \pmod{32}.
\end{align*}

Applying Lemmas~\ref{Syracuse:lemma:1} and~\ref{Syracuse:lemma:2} with $k=1$, 
we see that every positive integer (except integers $m \in \mcr_0$) in the congruence classes $13\pmod{32}$ and $5\pmod{64}$ 
has permutation pattern $(3,2,1)$.   
Applying Lemmas~\ref{Syracuse:lemma:3} and~\ref{Syracuse:lemma:4} with $k=1$, 
we see that  
every positive integer in the congruence classes $29\pmod{32}$ and $37\pmod{64}$ 
has permutation pattern $(3,1,2)$.   
Thus, the set of positive integers $m \notin \mcr_0$  in   
\[
5\pmod{8} \setminus  21\pmod{32} 
= 5\pmod{8} \setminus  r_2\pmod{2\cdot 4^2}  
\]
is partitioned into two sets, one with permutation pattern $(3,2,1)$ 
and the other with permutation pattern $(3,1,2)$.
Each of these sets has density
\[
\sum_{i=4}^5 \frac{1}{2^i} = \frac{1}{8}\left(1 -  \frac{1}{4} \right).
\]

The congruence class $21 \mod{32} = r_2 \pmod{32}$ is the disjoint union 
of the following five   congruence classes  (in the green boxes):
\begin{align*}
53\pmod{128} &  \\
21\pmod{256}  & \\
117\pmod{128} & \\
149 \pmod{256} & \\
85\pmod{128} & = r_3 \pmod{128}.
\end{align*}

Applying Lemmas~\ref{Syracuse:lemma:1} and~\ref{Syracuse:lemma:2} with $k=2$, 
we see that every positive integer  $m \notin \mcr_0$ 
in the congruence classes $53\pmod{128}$ and $21\pmod{256}$ 
has permutation pattern $(3,2,1)$.  
Applying Lemmas~\ref{Syracuse:lemma:3} and~\ref{Syracuse:lemma:4} with $k=2$, 
we see that 
every positive integer in the congruence classes $117\pmod{128}$ and $149 \pmod{256}$ 
has permutation pattern $(3,1,2)$.   
Thus, the set of positive integers $m \notin \mcr_0$  in 
\[
21\pmod{32} \setminus   85\pmod{128} 
= r_2\pmod{2\cdot 4^2} \setminus  r_3 \pmod{2\cdot 4^3}  
\]
is partitioned into two 
sets, one with permutation pattern $(3,2,1)$ and the other with permutation pattern $(3,1,2)$.
Each of these sets has density
\[
\sum_{i=6}^7 \frac{1}{2^i}  = \frac{1}{8}\left(  \frac{1}{4} -  \frac{1}{4^2} \right).
\]
It follows that the set of positive integers $m \notin \mcr_0$  in 
\[
5\pmod{8} \setminus  85\pmod{128}  
=  5\pmod{8}   \setminus  r_3 \pmod{2\cdot 4^3}  
\]
is partitioned into two sets, one with permutation pattern $(3,2,1)$ 
and the other with permutation pattern $(3,1,2)$.
Each of these sets has density
\[
\frac{1}{8}\left( 1-  \frac{1}{4^2} \right).
 \]

Let $k \geq 3$ and assume  that the set of  positive integers $m \notin \mcr_0$ in 
\[
5\pmod{8} \setminus r_{k} \pmod{2\cdot 4^{k}} 
\]
is partitioned into two sets, one with permutation pattern $(3,2,1)$ 
and the other with permutation pattern $(3,1,2)$, and that each of these sets has density
\[
\frac{1}{8}\left(1 -  \frac{1}{4^{k-1}} \right).
 \]

The following diagram displays the partition of the ``red'' congruence class 
$r_k \pmod{2\cdot 4^k}$ into three ``blue'' congruence classes 
modulo $2\cdot 4^{k+1}$ and two ``blue''  congruence classes 
modulo $4^{k+2}$.

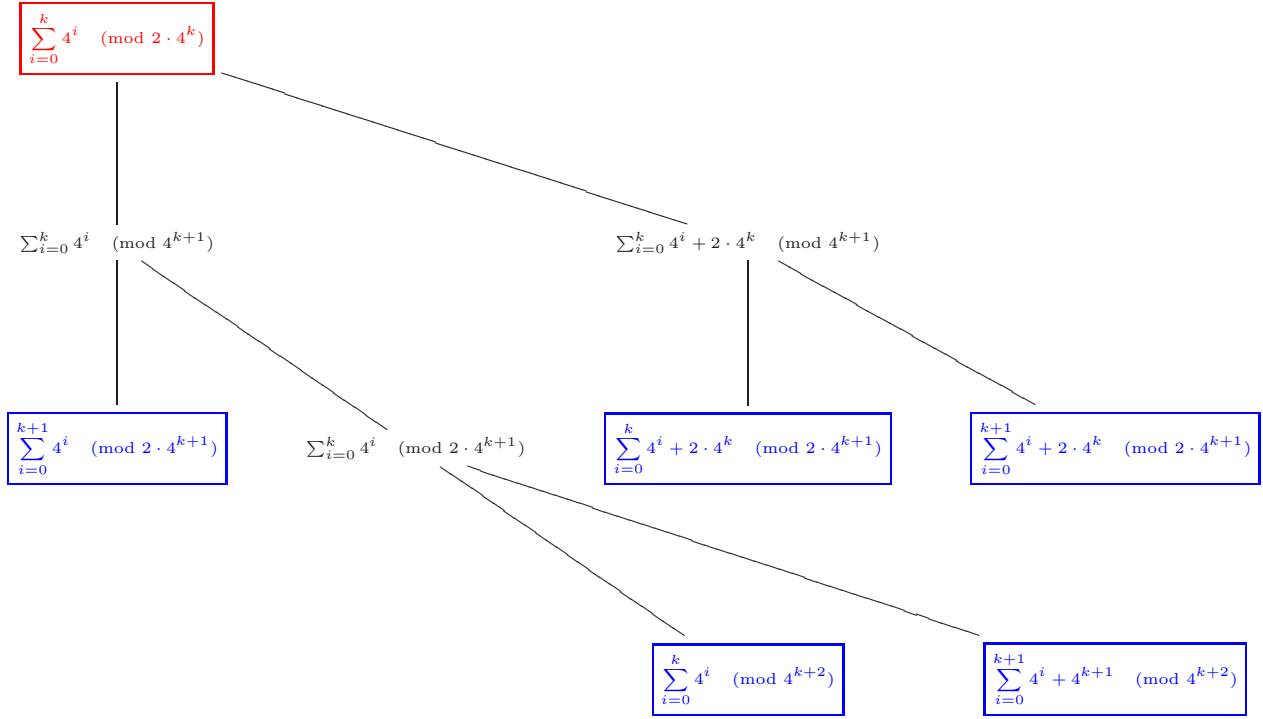
\tiny{
\[
\xymatrix{
 \color{red}\boxed{\sum_{i=0}^k 4^i \pmod{2\cdot 4^k}} \ar@{-}[dd]   \ar@{-}[ddrr]  & & &    \\
 & & &    \\
\sum_{i=0}^k 4^i \pmod{ 4^{k+1}}  \ar@{-}[dd]   \ar@{-}[ddr]    & & 
 \sum_{i=0}^k 4^i + 2 \cdot 4^k\pmod{ 4^{k+1}}  \ar@{-}[dd]   \ar@{-}[ddr] &  \\
& & &    \\
 \color{blue}\boxed{\sum_{i=0}^{k+1} 4^i \pmod{2\cdot 4^{k+1}} }
&  \sum_{i=0}^{k} 4^i \pmod{2\cdot 4^{k+1}}  \ar@{-}[ddr]   \ar@{-}[ddrr]  
  &  \color{blue}\boxed{  \sum_{i=0}^k 4^i + 2 \cdot 4^k\pmod{  2 \cdot  4^{k+1}}  }
  & \color{blue} \boxed{ \sum_{i=0}^{k+1} 4^i + 2 \cdot 4^k\pmod{  2 \cdot 4^{k+1}} } \\
  & & &    \\
 &   &  \color{blue}\boxed{  \sum_{i=0}^{k} 4^i \pmod{4^{k+2}} }  & \color{blue} \boxed{ \sum_{i=0}^{k+1} 4^i  + 4^{k+1}\pmod{4^{k+2}}}     \\
  }
  \]
 }
\normalsize{}

The congruence class $ r_k \pmod{2\cdot 4^k}$ is the disjoint union of the following five congruence classes:
\begin{align*}
\sum_{i=0}^k 4^i + 2 \cdot 4^k & \pmod{  2 \cdot  4^{k+1}} \\
\sum_{i=0}^k 4^i  &\pmod{    4^{k+2}} \\
\sum_{i=0}^{k+1} 4^i + 2 \cdot 4^k & \pmod{  2 \cdot  4^{k+1}} \\
\sum_{i=0}^{k+1} 4^i +  4^{k+1} & \pmod{   4^{k+2}} \\
\sum_{i=0}^{k+1} 4^i  & \pmod{  2 \cdot  4^{k+1}}.
\end{align*}

Applying Lemmas~\ref{Syracuse:lemma:1} and~\ref{Syracuse:lemma:2}, 
we see that every positive integer in the congruence classes 
$\sum_{i=0}^k 4^i + 2 \cdot 4^k  \pmod{  2 \cdot  4^{k+1}}$ and 
$\sum_{i=0}^k 4^i  \pmod{    4^{k+2}}$   
has permutation pattern $(3,2,1)$.  
Applying Lemmas~\ref{Syracuse:lemma:3} and~\ref{Syracuse:lemma:4}, 
we see that 
every positive integer in the congruence classes 
$\sum_{i=0}^{k+1} 4^i + 2 \cdot 4^k  \pmod{  2 \cdot  4^{k+1}}$ and 
$\sum_{i=0}^{k+1} 4^i +  4^{k+1}  \pmod{   4^{k+2}}$ 
has permutation pattern $(3,1,2)$.    
Thus, the positive integers in the set 
\[
r_k \pmod{2\cdot 4^k} \setminus r_{k+1} \pmod{2\cdot 4^{k+1}} 
\]
 are partitioned into two 
sets, one with permutation pattern $(3,2,1)$ and the other with permutation pattern $(3,1,2)$.
Each of these sets has density
\[
 \frac{1}{4^{k+1}} + \frac{2}{4^{k+2}} =  \frac{1}{8} \left( \frac{1}{4^{k-1}} - \frac{1}{4^k}\right).
\]

Thus, the   positive integers in  
\begin{align*}
5 & \pmod{8} \setminus  r_{k+1} \pmod{2\cdot 4^{k+1}}  \\
& = \left( 5\pmod{8} \setminus r_k \pmod{2\cdot 4^k} \right) \bigcup 
 \left( r_k \pmod{2\cdot 4^k} \setminus  r_{k+1} \pmod{2\cdot 4^{k+1}}  \right)
\end{align*} 
are partitioned into two sets, one with permutation pattern $(3,2,1)$ 
and the other with permutation pattern $(3,1,2)$.
Each of these sets has density
\[
  \frac{1}{8} \left( 1 - \frac{1}{4^{k-1}}\right) + \frac{1}{8} \left( \frac{1}{4^{k-1}} - \frac{1}{4^k}\right)
 =  \frac{1}{8} \left( 1 - \frac{1}{4^k}\right).
\] 
This completes the induction.  

For every $\varepsilon > 0$ there is an integer $k$ such that $1/(8\cdot 4^k) < \varepsilon$
and so both sets $\mcr_1$ and $\mcr_2$ 
contain subsets of density greater than $1/8 - \varepsilon$.  
The set $\mcr_0$ has density $0$ and the congruence class $5\pmod{8}$ 
has density $1/4$ with respect to $\Omega$. 
Applying Lemma~\ref{Syracuse:lemma:densityW} with $X = \Omega$, 
$W = \{ m>1: m \equiv 5\pmod{8}\}$, $t=2$, and $\alpha_1 = \alpha_2 = 1/8$ 
to the partition $5\pmod{8} = \mcr_0 \cup \mcr_1 \cup \mcr_2$, 
we see that $\mcr_1$ has density $1/8$ and $\mcr_2$ has density $1/8$. 
This completes the proof. 
\end{proof}


\section{Some impossible permutation patterns for quadruples} 
By Theorem~\ref{Syracuse:theorem:triples},
every triple permutation pattern is the permutation pattern of triples 
$\left( m, S(m), S^2(m) \right)$ 
of the iterated Syracuse function for a  set of integers  $m$ of positive density.
The story for quadruple permutation patterns is different.  
In this section we prove that there are quadruple permutation patterns 
that never occur as permutation patterns of quadruples 
$\left( m, S(m), S^2(m), S^3(m) \right)$ 
of the iterated Syracuse function. 

\bt                                    \label{Syracuse:theorem:quad-123}
Consider   quadruples 
\[
\left( m, S(m), S^2(m), S^3(m) \right)
\]
such that $S^i(m) \neq S^j(m)$ for $i \neq j$ and 
the triple $\left(m, S(m),  S^2(m) \right)$ has permutation pattern $(1,2,3)$, that is, 
\[
m< S(m) < S^2(m).
\]
For these quadruples there are four possible quadruple permutation patterns:
\[
(1,2,3,4), \qquad (1,2,4,3),    \qquad (2,3,4,1), \qquad (1,3,4,2).
\]
The density of each of these permutation patterns is as follows:
\begin{center}
\begin{tabular}{c|c}
permutation pattern $\sigma \in \Sigma_4$ & permutation density $d_4(\sigma)$ \\     \hline 
$(1,2,3,4)$ &  $1/8$ \\ 
$(1,2,4,3)$ & $1/16$ \\
$(2,3,4,1)$ & $1/16$ \\
$(1,3,4,2)$ &$0$  \\
\end{tabular}
\end{center}
Moreover, the permutation pattern $(1,3,4,2)$ never occurs.
\et

Note that $1/4 = 1/8 + 1/16 + 1/16$ is the Syracuse function permutation pattern density 
of the triple $(1,2,3)$.

\begin{proof}
By Theorems~\ref{Syracuse:theorem:triples-not-5} and~\ref{Syracuse:theorem:triples-5}, 
we have $m < S(m) < S^2(m)$ if and only if 
\[
m \equiv 7 \pmod{8}.
\]
Then $m = 7 + 8x$ for some nonnegative integer $x$ and  
\[
m = 7 + 8x < S(m) = 11+12x < S^2(m) = 17+18x.  
\]
It follows that 
\[
 S^3(m) = \frac{26 + 27x}{2^e}
\]
for some nonnegative integer $e$. 

The congruence class $7 \pmod{8}$ is the disjoint union of the congruence classes $15 \pmod{16}$, $7 \pmod{32}$, and $23 \pmod{32}$. 

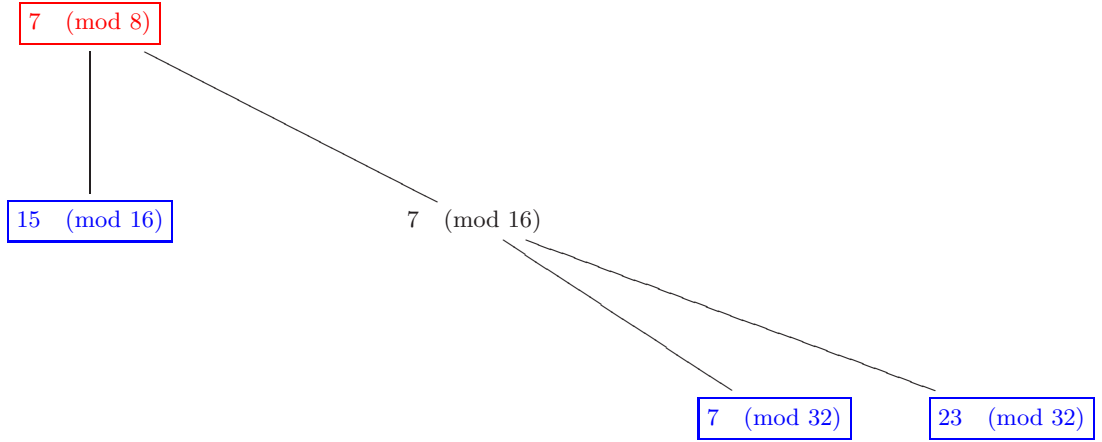
\small{
\[
\xymatrix{
\color{red}\boxed{ 7 \pmod {8} } \ar@{-}[dd]   \ar@{-}[ddrrr]  & & & & &  &  \\
 & & & & &  &   \\
 \color{blue}\boxed{ 15\pmod{16} }   &  & &  7 \pmod{16}  \ar@{-}[ddrr]   \ar@{-}[ddrrr] & & & \\
& & & & & &   \\
& & &  & &    \color{blue}\boxed{7 \pmod{32}}  & \color{blue}\boxed{23 \pmod{32}} \\
}
\]
}

If $m \equiv 15 \pmod{16}$, then the integers $x$ and  $26 + 27x$ are odd and so $e=0$.  
We have 
\[
m < S(m) < S^2(m) = 17+18x <  26 + 27x = S^3(x)
\]
and so the quadruple $(m, S(m), S^2(m), S^3(m))$ has permutation pattern $(1,2,3,4)$. 

If $m \equiv 7 \pmod{32}$, then  $x = 4y$ and  $e=1$.  
We obtain 
\[
S^3(m) = \frac{26 + 27x}{2}  = 13 + 54y = 13 + \left(\frac{27}{2}\right) x.
\]
The inequality 
\[
7 + 8x < 11+12x < 13 + \left(\frac{27}{2}\right) x <  17+18x
\]
implies 
\[
m < S(m) < S^3(m) < S^2(m). 
\]
The quadruple $(m, S(m), S^2(m), S^3(m))$ has permutation pattern $(1,2,4,3)$.

If $m \equiv 23 \pmod{32}$, then  $x = 2+4y$ and $e\geq 2$.   
\begin{align*}
S^3(m)&  = \frac{26 + 27x}{2^e}   = \frac{80 + 27\cdot 4y}{2^e} 
 = \frac{20 + 27 y}{2^{e-2}} \\ 
 & \leq 20 + 27 y = 20 + 27\left(\frac{x-2}{4}\right) \\
& = \frac{13}{2} + \left(\frac{27}{4}\right)x < 7+8x = m.
\end{align*}
The quadruple $(m, S(m), S^2(m), S^3(m))$ has permutation pattern 
$(2,3,4,1)$.

We see that the permutation pattern $(1,3,4,2)$ never occurs, and that the 
permutation patterns $(1,2,3,4)$, $(1,2,4,3)$, and $(2,3,4,1)$ have 
permutation densities $1/8$, $1/16$, and $1/16$, respectively. 
This completes the proof. 
\end{proof}

\bt                                        \label{Syracuse:theorem:quad-132}
Consider  quadruples 
\[
\left( m, S(m), S^2(m), S^3(m) \right)
\]
such that $S^i(m) \neq S^j(m)$ for $i \neq j$ and 
the triple $\left(m, S(m),  S^2(m) \right)$ has permutation pattern $(1,3,2)$, that is,
\[
m < S^2(m) < S(m).
\]  
For these quadruples there are four possible quadruple permutation patterns:
\[
(1,3,2,4), \qquad (2,4,3,1), \qquad (1,4,2,3),    \qquad (1,4,3,2) .
\]
The density of each of these permutation patterns is as follows:
\begin{center}
\begin{tabular}{c|c}
permutation pattern $\sigma \in \Sigma_4$ & permutation density $d_4(\sigma)$ \\   \hline 
$( 1,3,2,4)$ &  $1/16 $  \\ 
$( 2,4,3,1) $ &  $ 1/16$  \\
$( 1,4,2,3)$ & $ 0$  \\
$( 1,4,3,2)$ & $ 0$  
\end{tabular}
\end{center}
Moreover, the permutation patterns $( 1,4,3,2)$ and $( 1,4,2,3)$ never occur.
\et

Note that $1/8 = 1/16 + 1/16$ is the Syracuse function permutation pattern density 
of the triple $(1,3, 2)$.  

\begin{proof}
By Theorems~\ref{Syracuse:theorem:triples-not-5} and~\ref{Syracuse:theorem:triples-5},
we have $m < S^2(m) < S(m)$ if and only if 
\[
m \equiv 11 \pmod{16}. 
\]
Then $m = 11 + 16x$ for some nonnegative integer $x$ and  
\[
m = 11 + 16x < S^2(m) = 13 + 18x <    S(m) = 17+24x.
\]
It follows that 
\[
 S^3(m) = \frac{20 + 27x}{2^e}
\]
for some nonnegative integer $e$.

The congruence class $11 \pmod{16}$ is the disjoint union of the congruence 
classes $11 \pmod{32}$  and $27 \pmod{32}$. 

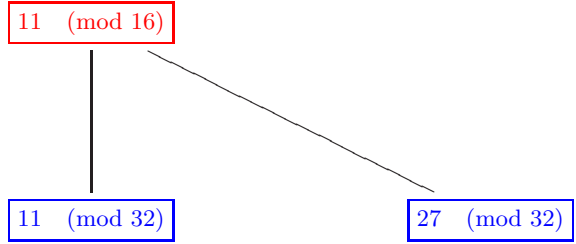
\small{
\[
\xymatrix{
\color{red}\boxed{ 11 \pmod {16} } \ar@{-}[dd]   \ar@{-}[ddrrr]  & & & & &  &  \\
 & & & & &  &   \\
 \color{blue}\boxed{ 11\pmod{32} }   &  & &   \color{blue}\boxed{27 \pmod{32} }& & & 
}
\]
}

If $m \equiv 27 \pmod{32}$, then $x$ and $20 + 27x$ are odd and so  $e=0$. 
We have 
\[
 S(m) = 17 + 24x <  20 + 27x = S^3(m) 
\]
and the quadruple $(m, S(m), S^2(m), S^3(m))$ has permutation pattern $(1,3,2,4)$. 

If  $m \equiv 11 \pmod{32}$, then $x$ is even and $e \geq 1$.  We obtain  
\begin{align*}
 S^3(m)  & = \frac{20 + 27x}{2^e} \leq 10 + \left( \frac{27}{2} \right)  x  < 11 + 16x = m. 
\end{align*}
The quadruple $(m, S(m), S^2(m), S^3(m))$ has permutation pattern $(2,4,3,1)$. 

We see that the permutation patterns $( 1,4,3,2)$ and $( 1,4,2,3)$ never occur 
and that each of the 
permutation patterns $( 1,3,2,4)$ and $( 2,4,3,1) $ has density $1/16$. 
This completes the proof. 
\end{proof}

\bt                                   \label{Syracuse:theorem:quad-213}
Consider  quadruples 
\[
\left( m, S(m), S^2(m), S^3(m) \right)
\]
such that $S^i(m) \neq S^j(m)$ for $i \neq j$ and 
the triple $\left(m, S(m),  S^2(m) \right)$ has permutation pattern $(2,1,3)$, that is,
\[
S(m) < m < S^2(m).
\]  
For these quadruples there are four possible quadruple permutation patterns:
\[
(2,1,3,4), \qquad (2,1,4,3),    \qquad (3,1,4,2), \qquad (3,2,4,1) .
\]
The density of each of these permutation patterns is as follows:
\begin{center}
\begin{tabular}{c|c}
permutation pattern $\sigma \in \Sigma_4$ & permutation density $d_4(\sigma)$ \\    \hline 
$(2,1,3,4)$ &  $1/16$  \\ 
$(3,1,4,2)$ & $1/32$  \\
$(3,2,4,1) $ &  $1/32$  \\
$(2,1,4,3))$ & $0$  \\
\end{tabular}
\end{center}
Moreover, the permutation pattern $(2,1,4,3 )$  never occurs.
\et

Note that $1/8 = 1/16 + 1/32 + 1/32$ is the Syracuse function permutation pattern density 
of the triple $(2,1,3)$.  

\begin{proof}
By Theorems~\ref{Syracuse:theorem:triples-not-5} and~\ref{Syracuse:theorem:triples-5},
we have $S(m) < m < S^2(m)$ if and only if 
\[
m \equiv 9 \pmod{16}. 
\]
Then $m = 9 + 16x$ for some nonnegative integer $x$ and  
\[
S(m) = 7+12x <  m = 9 + 16x < S^2(m) = 11+18x.
\]
It follows that 
\[
 S^3(m) = \frac{17 + 27x}{2^e}
\]
for some nonnegative integer $e$. 

The congruence class $9 \pmod{16}$ is the disjoint union of the congruence classes 
$9 \pmod{32}$, $25 \pmod{64}$, and $57 \pmod{64}$. 

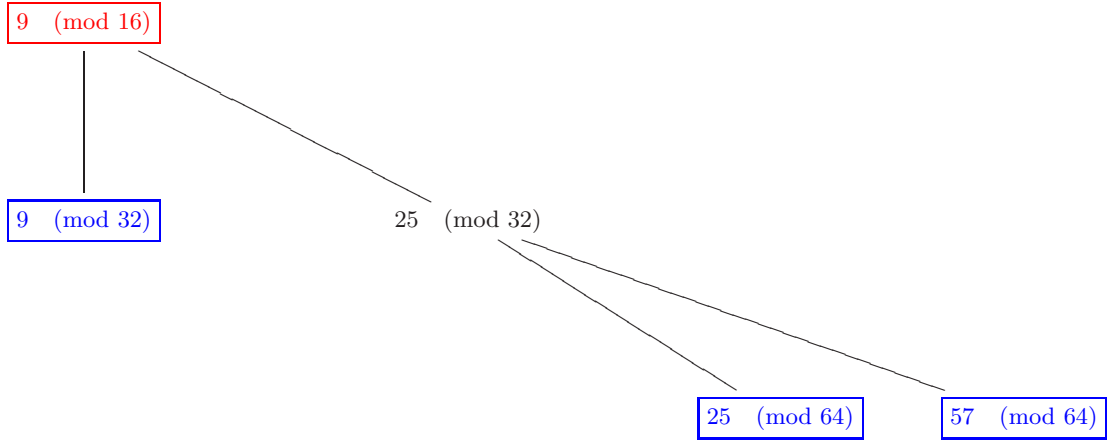
\small{
\[
\xymatrix{
\color{red}\boxed{ 9 \pmod {16} } \ar@{-}[dd]   \ar@{-}[ddrrr]  & & & & &  &  \\
 & & & & &  &   \\
 \color{blue}\boxed{9\pmod{32} }   &  & &  25 \pmod{32}  \ar@{-}[ddrr]   \ar@{-}[ddrrr] & & & \\
& & & & & &   \\
& & &  & &    \color{blue}\boxed{25 \pmod{64}}  & \color{blue}\boxed{57 \pmod{64}} \\
}
\]
}

If $m \equiv 9 \pmod{32}$, then $x$ is even,  $17 + 27x$ is odd, and so $e=0$.  
We have 
\[
S^2(m) = 11+18x < S^3(m) = 17 + 27x 
\]
and the quadruple $(m, S(m), S^2(m), S^3(m))$ has permutation pattern $(2,1,3,4)$. 

If $m \equiv 57 \pmod{64}$, then $x = 3+4y$ and $e=1$.  We obtain 
\begin{align*}
S^3(m) & = \frac{17 + 27(3+4y)}{2} = 49+54y \\
& = 49+54 \left( \frac{x-3}{4}  \right) 
= \frac{17}{2}+ \left( \frac{27}{2}\right) x. 
\end{align*}
The inequality
\[
7+12x <\frac{17}{2}+ \left( \frac{27}{2}\right) x < 9+16x
\]
implies 
\[
S(m) < S^3(m) < m < S^2(m).
\]
The quadruple $(m, S(m), S^2(m), S^3(m))$ has permutation pattern $(3,1,4,2)$.

If $m\equiv 25 \pmod{64}$,  then $x = 1+4y$ and  $e \geq 2$. 
We obtain  
\begin{align*}
S^3(m) 
& = \frac{17 + 27(1+4y)}{2^e} = \frac{11 + 27y}{2^{e-2}} \\
& \leq 11 + 27y = 11 + 27\left(\frac{x-1}{4} \right)  = \frac{17}{4}+ \left( \frac{27}{4}\right) x \\
& < 7+12 x = S(m). 
\end{align*}
The quadruple $(m, S(m), S^2(m), S^3(m))$ has permutation pattern $(3,2,4,1)$. 

We see that the permutation pattern $(2,1,4,3 )$  never occurs, and that the 
permutation patterns $(2,1,3,4)$, $(3,1,4,2)$, and $(3,2,4,1)$ have densities 
$1/16$, $1/32$, and $1/32$, respectively. 
This completes the proof. 
\end{proof}


\section{Dropping time}
An odd integer $m > 1$ has \emph{dropping time}\index{dropping time} $D(m) = k$ 
if $k$ is the smallest positive integer such that $S^k(m) < m$, 
and dropping time $D(m) = \infty$ if $S^k(m) > m$ for all positive integers $k$. 
The Collatz conjecture is equivalent to the statement that $D(m) < \infty$ for 
all odd integers $m > 1$. 
We have $D(m) > k$ if and only if $S^i(m) \neq 1$ for all $i \in \{1,2,\ldots, k\}$ 
and the permutation pattern of the $(k+1)$-tuple 
$\left(m, S(m), S^2(m), \ldots, S^k(m) \right)$ is not of the form 
$(1,a_2,a_3,\ldots, a_{k+1})$, where $(a_2,a_3,\ldots, a_{k+1})$ is any permutation of 
$(2,3,\ldots, k+1)$.

\bt
Let $N_k(x)$ count the number of odd integers $m \leq x$ such that $D(m) \leq k$. 
Then 
\begin{align*}
\overline{N}_1 = \lim_{x\rightarrow \infty} \frac{N_1(x)}{x/2} & = \frac{1}{2} \\
\overline{N}_2 =\lim_{x\rightarrow \infty} \frac{N_2(x)}{x/2}  & = \frac{5}{8} \\
\overline{N}_3 = \lim_{x\rightarrow \infty} \frac{N_3(x)}{x/2}  & = \frac{3}{4} \\
\end{align*}
\et

\begin{proof}
Permutation pattern densities give explicit values for $N_k(x)$ 
for $k=1, 2$, and $3$.
An odd integer $m > 1$ has dropping time $D(m) = 1$ if and only if $S(m) < m$
if and only if the pair $(m,S(m))$ has permutation pattern $(2,1)$, and so $\overline{N}_1$ 
is the density of the permutation pattern $(2,1)$.  
We have $\overline{N}_1 = 1/2$ by Theorem~\ref{Syracuse:theorem:pairs}. 
 
 The odd integer $m > 1$ has dropping time $D(m) \leq 2$ 
if and only if the permutation pattern of the  triple $(m,S(m),S^2(m))$ is not 
 $(1,2,3)$ or $(1,3,2)$.  
By Theorem~\ref{Syracuse:theorem:triples}, the permutation pattern $(1,2,3)$ 
has density $1/4$ and   the permutation pattern $(1,3,2)$ 
has density $1/8$.  
Therefore, $\overline{N}_2  = 1 - 1/4 - 1/8 = 5/8$. 

 The odd integer $m > 1$ has dropping time $D(m) \leq 3$ 
if and only if the quadruple $(m,S(m),S^2(m), S^3(m))$ has permutation pattern not 
equal to $(1, a_2,a_3,a_4)$, where $(a_2,a_3,a_4)$ is any 
of the six permutations of $2,3$, and $4$. 
By Theorems~\ref{Syracuse:theorem:quad-123} and~\ref{Syracuse:theorem:quad-132}, 
the sum of the densities of these six permutation patterns is 1/4 
and so $\overline{N}_3  = 1 - 1/4 = 3/4$. 
\end{proof}

The dropping time function is the  Syracuse function analogue of 
the stopping time function of Riho Terras~\cite{terr76}.  
It would be of interest to prove (similar to results of Terras) that the limit 
\[
\overline{N}_k = \lim_{x\rightarrow \infty} \frac{N_k(x)}{x/2} 
\]
exists for all positive integers $k$ and that 
\[
\lim_{k \rightarrow \infty} N_k = 1. 
\]

\textbf{Acknowledgement.} 
I thank Yosef Berman (CUNY Graduate Center) for computer calculations 
of permutation patterns for triples and quadruples of the iterated Syracuse 
function $S(m)$ that motivated this paper.

\end{document}